\documentclass{amsart}
\usepackage{amssymb,latexsym}
\usepackage{amsmath}
\usepackage{graphicx}

\newtheorem{theorem}{Theorem}[section]

\newtheorem{corollary}[theorem]{Corollary}
\newtheorem*{cor}{Corollary 5.5'}
\newtheorem{remark}[theorem]{Remark}
\newtheorem{lemma}[theorem]{Lemma}

\newtheorem{observation}[theorem]{Observation}

\newtheorem*{thm}{Theorem(Wu, [5])}

\newtheorem{defin}[theorem]{Definition}

\newcommand{\bbb}{\mbox{$\beta$}}
\newcommand{\aaa}{\mbox{$\alpha$}}
\newcommand{\eee}{\mbox{$\epsilon$}}
\newcommand{\Rrr}{\mbox{$\mathbb{R}$}}
\newcommand{\bdd}{\mbox{$\partial$}}

\begin{document}


\title{Compressing thin spheres in the complement of a link }
\author{Maggy Tomova}
\address{Department of Mathematics, University of California, Santa
Barbara, CA 93117} \email{maggy@math.ucsb.edu}
\thanks{Research supported in part by an NSF grant}
\thanks{Version 2}

\begin{abstract}
Let $L$ be a link in $S^3$ that is in thin position but not in
bridge position and let $P$ be a thin level sphere. We generalize a result of Wu by giving a bound on
the number of disjoint irreducible compressing disks $P$ can have,
including identifying thin spheres with unique compressing disks. We
also give conditions under which $P$ must be incompressible
on a particular side or be weakly
incompressible. If $P$ is strongly compressible we describe how a
pair of compressing disks must lie
relative to the link.
\end{abstract}
\maketitle
\section{Introduction}

For $L$ a link in $S^3$, the natural height function $p:(S^3-poles)
\rightarrow \Rrr$ gives rise to level spheres which correspond to
meridional planar surfaces in the complement of $L$. The link
is said to be in thin position if the number of intersection
points between the level spheres and the link has been minimized;
the precise definition is given later. A level sphere is called {\em thin}
if the highest critical point for $L$ below it is a maximum and the
lowest critical point above it is a minimum. The link is said to be in
bridge position if it does not have any thin spheres. In \cite{Thom}, Thompson
has shown that if a knot is in thin, but not in bridge position,
maximally compressing a thin level sphere results in a non-trivial
incompressible meridional planar surface.  Wu \cite{Wu} has shown
that a thin sphere of
minimum width is itself incompressible in the link complement. A
natural question to ask is what can be said about the compressing
disks for other thin spheres in $S^3-L$ given information
about the spheres' widths. In this paper we give some sufficient conditions for a thin sphere to be
incompressible on a particular side, to be weakly incompressible, or
to have a unique compressing disk. If $P$ is strongly compressible we describe
how a pair of
compressing disks must be positioned relative to the link. We also
prove a bound on the number of disjoint ``simple" compressing disks a
thin sphere can have based on its width, thus generalizing Wu's result.
In the last section we give an additional restriction on the
compressibility of a thin sphere in the case where the link is prime or is a knot.

\section{Thin position}

We start with a quick review of thin position, a concept
originally due to Gabai \cite{Ga}. A detailed overview can be
found in \cite{Schar1}. Let $p: S^3 \to \mathbb{R}$ be the
standard height function and let $L$ be a link in $S^3$ such that
$p$ restricts to a Morse function on $L$. If $t$ is a regular
value of $p|L$, $p^{-1}(t)$ is called a level sphere with width
$w(p^{-1}(t))=|L\cap p^{-1}(t)|$. If $c_{0}<c_{1}<...<c_{n}$ are all the
critical values of $p|L$, choose regular values
$r_{1},r_{2},...,r_{n}$ such that $c_{i-1}<r_{i}<c_{i}$. Then the
{\em width of $L$ with respect to $p$} is defined by $w(L,p)=\sum
w(p^{-1}(r_{i}))$. The {\em width} of $L$, $w(L)$ is the minimum of $w(L',p)$
over all $L'$ isotopic to $L$. We say that $L$ is in thin position
if $w(L,p)=w(L)$.
\begin{remark}\label{rmk:sliding that thins}
\textnormal {Isotoping $L$ so as to slide a minimum above a maximum
decreases the width by 4. Sliding a minimum below a maximum
increases the width by 4, while sliding a minimum past a minimum or
a maximum past a maximum has no effect on the width.}
\end{remark}
If $\gamma \subseteq L$ is a 1-manifold (not necessarily
connected) and $P$ is a level sphere, then we call $P$ a
\textit{thin} (resp. \textit{thick}) \textit{level sphere for
$\gamma$} if the lowest critical value of $p|\gamma$ above $P$ is
a minimum (resp. maximum) and the highest critical value below $P$
is a maximum (resp. minimum). To avoid having to consider special
cases later, we will define $P$ to be thin for $\gamma$ also when
$\gamma\cap P=\emptyset$ . If $\gamma=L$ we call $P$ thin or
thick without specifying the 1-manifold. We say that $L$ is in bridge
position if there are no thin spheres.

Since $L$ is in general position with respect to $p$ it is
disjoint from both the minimum (south pole) and maximum (north
pole) of $p$ on $S^{3}$.  The complement of these poles in $S^{3}$
is of course diffeomorphic to $S^{2} \times \Rrr$ and, using $(-1,
1) \cong \Rrr$, we can choose the diffeomorphism to preserve the
foliation of $S^{3} - poles$ by level spheres given by $p$.  The
width of $L$ could just as easily be computed via its
diffeomorphic image in $S^{2} \times \Rrr$.  So, with little risk
of substantive confusion, we will often regard $L$ as contained in
$S^{2} \times \Rrr$ and continue to use $p$ to denote the
projection on the second factor $p: S^{2} \times \Rrr \rightarrow
\Rrr.$

\section{Moving $L$ around in 3-space}

It will be useful to move parts of $L \subset S^{2} \times \Rrr$
vertically, that is without changing the projection of $L$ to
$S^{2}$, but altering only the height function $p$ on those parts.
Suppose, for example, $a < b$ are regular values for $p|L$.  Take
$\epsilon > 0$ so small that there are no critical values of $p|L$
in either of the intervals $[a, a+\epsilon]$ or $[b, b +
\epsilon]$.  Let $h: [a, b + \epsilon] \rightarrow [a, b +
\epsilon]$ be the automorphism that consists of the union of the
linear homeomorphisms $[a, a+ \epsilon] \rightarrow [a, b]$ and
$[a + \epsilon, b + \epsilon] \rightarrow [b, b + \epsilon].$

\begin{figure}
\begin{center} \includegraphics[scale=.5]{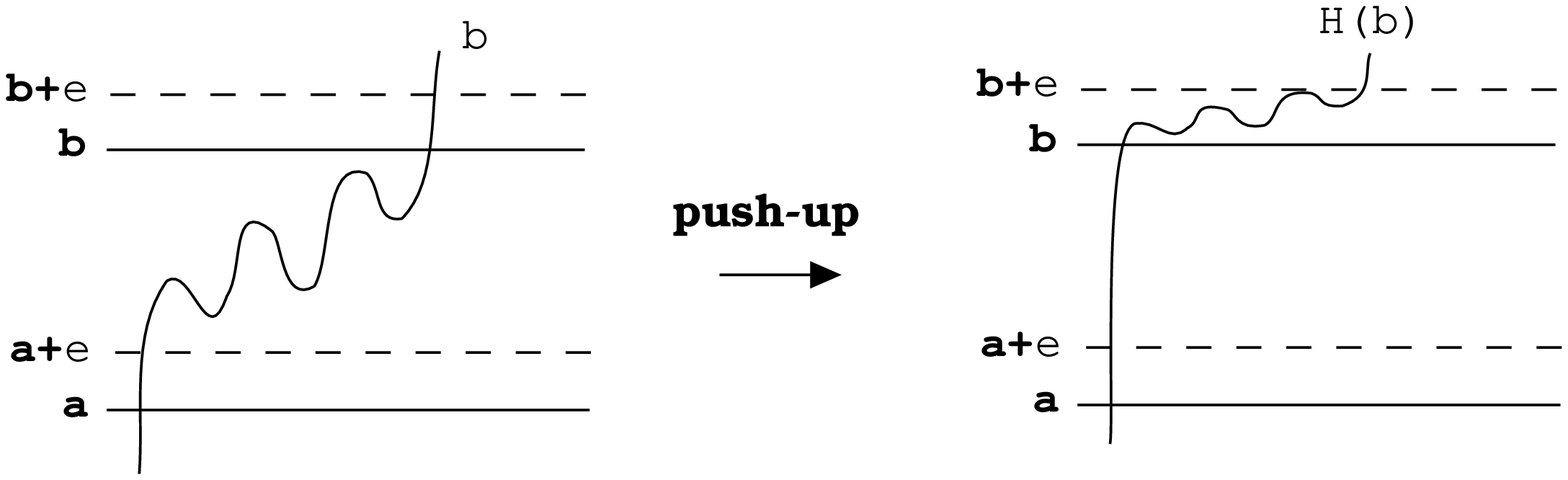}
 \end{center}
\caption{} \label{fig:pushup} \end{figure}

\begin{defin} Let $\beta$ be a collection of components of $L \cap
(S^{2} \times [a, b + \epsilon])$.  The push-up of $\beta$ past
$S^{2} \times \{ b \}$ is the image of $\beta$ under the
homeomorphism $H: S^{2} \times [a, b + \epsilon] \rightarrow S^{2}
\times [a, b + \epsilon]: (x, t) \mapsto (x, h(t))$. (Figure
\ref{fig:pushup})
\end{defin}

Notice that all critical points of $H(\beta)$ lie in $S^{2} \times
[b, b+\epsilon]$.  Since there is a linear isotopy from
$h$ to the identity, $\beta$ is properly isotopic to $H(\beta)$ in
$S^{2} \times [a, b+\epsilon]$.  This isotopy from $\beta$ to
$H(\beta)$ is called {\em pushing the critical points of $\beta$
above the sphere $p^{-1}(b)$}.  There is an obvious symmetric
isotopy that pushes the critical points of $\beta$ below the
sphere $p^{-1}(a)$.  These isotopies of $\beta$ only make sense as
isotopies of $L$ if they do not move $\beta$ across any
other part of $L$.  To that end, we now enhance somewhat a
fundamental construction due to Wu \cite{Wu}.

\begin{lemma} \label{lemma:Wu}
Let $P$ be a level sphere for $L \subset S^{3}$ and let $B$ be a
closed component of $S^{3} - P$.  Suppose $D$ is a compressing
disk for $P - L$ in $B - L$.  Then there is an isotopy of $L \cap
B$ rel $L \cap P$ in $B$ so that, after the isotopy, $p|L$ is
unchanged and $L \cap B$ is disjoint from $\bdd D \times \Rrr
\subset S^{2} \times \Rrr \subset S^{3}$.
\end{lemma}

\begin{remark}
\textnormal {An isotopy that, in the end, leaves $p|L$ unchanged is called
an $h$-isotopy in \cite{Wu}.  Adding the adjacent pole to $B \cap
(\bdd D \times \Rrr)$ creates a {\em vertical} disk in $B$ that is
a compressing disk for $P - L$ in $B - L$ and that has the same
boundary as $D$.}
\end{remark}

\begin{proof} The proof extends the argument of Lemma 2 of \cite{Wu}.
With no loss of generality we may assume that $P = p^{-1}(0)$, $B$
and $D$ lie above $P$, and $\bdd D$ is the equator of $P$,
dividing $P$ into two disks, $D_{w}$ and $D_{e}$.  Let $w$ and $e$ (the west
and east poles) denote the centers of $D_{w}$ and $D_{e}$
respectively.  Let $\gamma_{w}$ and $\gamma_{e}$ denote the unions of
the north pole with respectively the vertical arcs $\{ w \}
\times [0, \infty)$ and $\{ e \}
\times [0, \infty)$ which are contained in $S^{2} \times [0, \infty) \subset B$.
(Figure \ref{fig:Hisotopy} represents the situation one dimension
lower.) By general position, we may assume that $\gamma_{w}$ and
$\gamma_{e}$ are disjoint from $L$. By possibly piping points of $D
\cap \gamma_{e}$ over the pole we may assume, without moving $L$,
that $D$ is disjoint from $\gamma_{e}$.  In particular,  both
$D$ and $L \cap B$ lie inside $(S^{2} - \{ e \}) \times [0,
\infty) \cong \Rrr^{2} \times [0, \infty)$.  Choose the
parameterization $(S^{2} - \{ e \}) \cong \Rrr^{2}$ so that $\bdd
D$ is the unit circle and $\{ w \} = 0 \in \Rrr^2$.  (Then
$\gamma_{w}$ corresponds to $\{ 0 \} \times [0, \infty)$.)

We are now in a position to apply Wu's argument almost verbatim,
but adapted to our parameterization: Let $B_{w}$ and $B_{e}$ be the
balls in $B$ bounded by $D \cup D_{w}$ and $D \cup D_{e}$
respectively.  Let $\alpha = L \cap B_{w}$ and $\beta = L \cap B_{e}$.
Perform a level-preserving isotopy (fixed near $L \cap P \cong L
\cap (\Rrr^{2} \times \{ 0 \})$) that shrinks at every level
the radial distance from $0$ in $\Rrr^{2}$ so that, after the isotopy,
$B_{w} \subset D_{w} \times [0, \infty)$.  Next, without moving
$\beta$, shrink the ball $B_{w}$ so that it is very close to
$\Rrr^{2} \times \{ 0 \}$.  (This isotopy is not level-preserving
on $\alpha$.)  
\begin{figure}
\begin{center} \includegraphics[scale=.5]{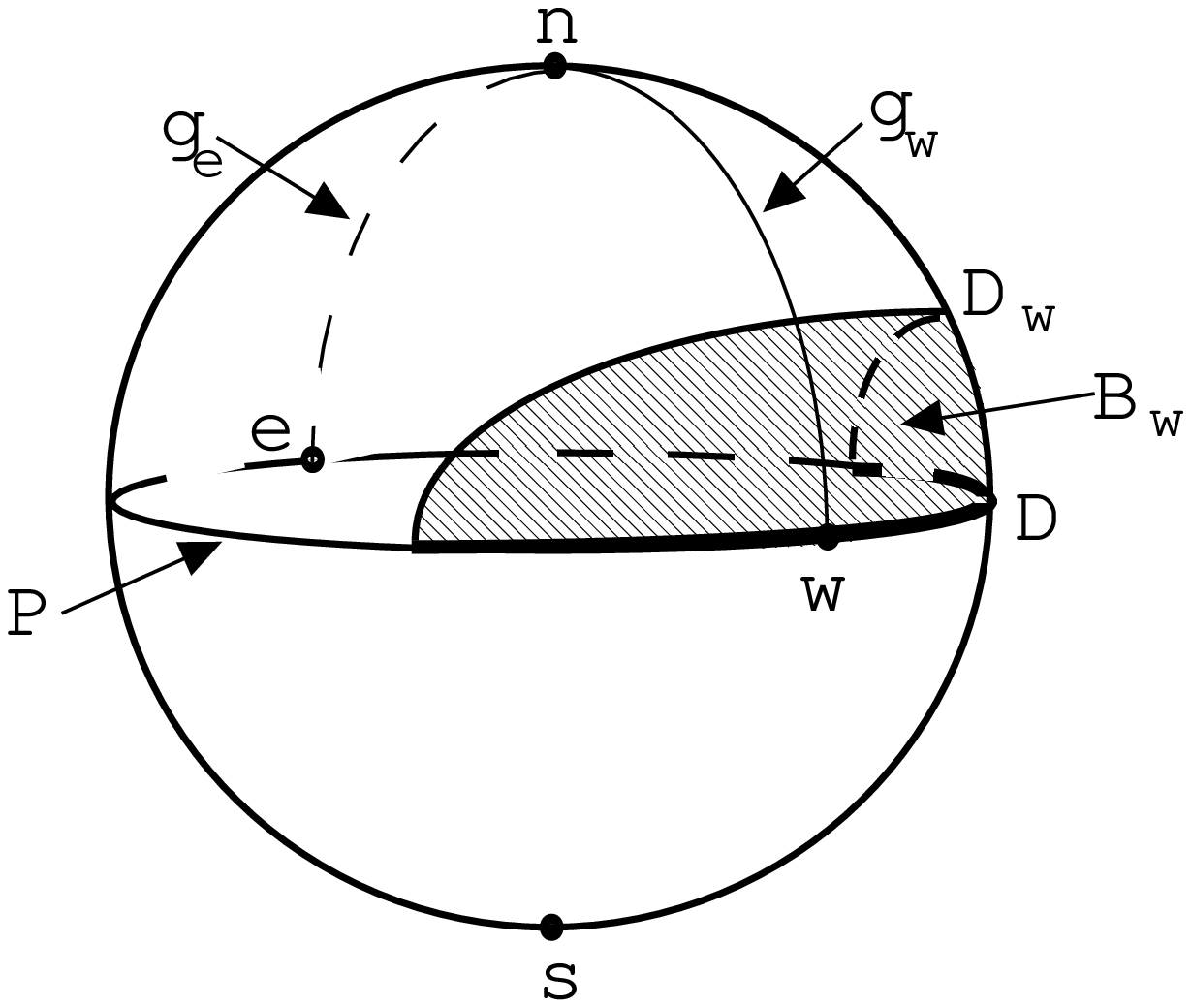} \end{center}
\caption{} \label{fig:Hisotopy}
\end{figure}

Then perform a level-preserving isotopy on $\beta$,
increasing at every level its radial distance from $0$ in
$\Rrr^{2}$, until $\beta$ lies entirely outside of $D_{w} \times
[0, \infty)$.  Finally, undo the isotopy that shrank $B_{w}$,
restoring $\alpha$ to its original height, but leaving it
entirely inside $D_{w} \times \Rrr$.
\end{proof}

It is straightforward to see that if a split link is placed in thin
position then there is a splitting sphere that is also a level sphere.
Thus little is lost in assuming henceforth that no link is split.
For the remainder of this paper we will assume that $L$ is an unsplit link
in $S^3$ in thin position and $P$ is a thin sphere compressible in
$S^3-L$ with compressing disk $D$. Without loss of generality we
assume $P=p^{-1}(0)$ and $D \subset p^{-1}[0,\infty)$. Then $D$
bounds two balls, $B_e$ and $B_w$ in $p^{-1}[0,\infty)$. Let
$\alpha=B_e \cap L$ and $\bbb=B_w \cap L$. We will be examining
the heights of the critical points for $\alpha$ and $\beta$ and
studying how frequently $\aaa$ and $\bbb$ intersect level spheres.
Both are unaffected by $h$-isotopies of $L$, so, following Lemma
\ref{lemma:Wu}, we may also assume that $\aaa$ and $\bbb$ are
separated by $\bdd D \times \Rrr$ and thus can be isotoped
vertically independently of one another.

\begin{lemma} \label{lemma:Thick levels don't overlap}

Let $T=p^{-1}(t)$ be a thick sphere for $\alpha$ and let
$S_+=p^{-1}(s_+)$ and $S_-=p^{-1}(s_-)$ be the first thin spheres
for $\alpha$ above and below $T$ respectively. Then, $T$ is a thin
sphere for $\beta$ and $\beta$ has no critical points in
$p^{-1}[s_-,s_+]$.

\end{lemma}

\begin{proof}
We first
show that $\beta$ does not have a minimum in the region
$p^{-1}(t,s_+)$. Isotope a small collar $\beta\cap
p^{-1}(t,t+\eee)$ upwards, pushing all critical points of
$\beta\cap p^{-1}(t,s_+)$ above $S^+$. Because $\alpha \cap p^{-1}(t,s_+)$
only has maxima, this isotopy slides all critical points of
$\beta\cap p^{-1}(t,s_+)$ only above maxima of $\alpha$. If any
critical point of $\beta$ is a minimum, this would decrease the
width of $L$ by Remark \ref{rmk:sliding that thins} contradicting
 $L$ being in thin position. Similarly, $\beta\cap p^{-1}(s_-,t)$
cannot have any maxima.

It follows that $\bbb \cap p^{-1}(s_-,s_+)$ has no
critical points at all. For suppose it has a minimum, necessarily
below T. Push $\beta \cap p^{-1}(s_-, s_-+\eee)$ upwards, above
$S_+$. Since $\beta$ has no maximum below $T$ and $\aaa$ has no
minimum above $T$ in $p^{-1}(s_-,s_+)$, no maximum is pushed up
past a minimum and the minimum of $\bbb \cap p^{-1}(s_-,t)$ is
moved past the maxima in $\aaa \cap p^{-1}(t,s_+)$ thinning $L$, a
contradiction. A similar argument shows that $\bbb \cap
p^{-1}(s_-,s_+)$ does not have any maxima and therefore $\bbb \cap
p^{-1}(s_-,s_+)$ is a product.

Finally, we show $T$ is a thin sphere for $\bbb$.  Suppose, to the
contrary, that the highest critical point of $\beta$ below $S_-$
is a minimum. Isotope $\beta$ so the minimum moves above $S_+$
making use of the product structure of $\beta$ between $S_-$ and
$S_+$. That would result in sliding a minimum up past the maxima
for $\alpha$, thus thinning $L$. Therefore the highest critical
point of $\beta$ below $T$ must be a maximum. Similarly, the lowest
critical point of $\beta$ above $T$, if it exists, must be a minimum proving the
lemma.

\end{proof}
\section{Alternating Levels}
Recall our assumption that $L$ is in thin position, $P=p^{-1}(0)$
is a compressible punctured sphere with compressing disk $D$ lying
above $P$, and $\aaa$ and $\beta$ are the strands of $L$ on each
side of $D$. In the following choose the labels $\aaa$ and $\bbb$
so that the maximum height of $\beta$ is greater than the maximum
height of $\alpha$. Let $A=p^{-1}(a)$ be the first thin sphere
above $\alpha$. If  $C=p^{-1}(c)$ is a thin sphere with $0<c \leq a$
we call $C$ an \textit{alternating sphere} for $D$ and $c$ an
\textit{alternating level} for $D$ if the minimum just above it and
the maximum just below it are on different sides of the disk $D$.

\begin{remark}\label{rmk:top is alternating}
\textnormal{It is always the case that $A$ is an alternating sphere for D.}
\end{remark}

\begin{remark}\label{rmk:product region}
\textnormal{Suppose $c'<c$ are two adjacent alternating levels. Then one
of $\alpha \cap p^{-1}[c',c]$ or $\beta \cap p^{-1}[c',c]$ is a
product. If $c$ is the lowest alternating level above $P$, then
$\alpha \cap p^{-1}[0,c]$ or $\beta \cap p^{-1}[0,c]$ is a
product.}
\end{remark}

\begin{lemma}\label{lemma:Top alternating sphere is thinnest}
For any level sphere $S=p^{-1}(s)$ with $s \in [0,a)$, $|A\cap
\beta|\leq|S\cap \beta|$ with equality if and only if $\beta \cap
p^{-1}(s,a)$ is a product. It follows that $|A\cap L|<|S\cap L|$.

\end{lemma}
\begin{proof}
 If $S$ is a level sphere
providing a counterexample to the claim, then at least one of the
two adjacent thin spheres is also a
counterexample, so without loss of generality we may assume that
$S$ is thin. Take $S=p^{-1}(s)$ to be the highest thin sphere
between $P$ and $A$ such that $|A\cap \beta|\geq|S\cap \beta|$ and
$\beta \cap p^{-1}(s,a)$ is not a product. Let $Q_1,..,Q_k$ be the
level spheres coming from critical values for $\alpha$ above $S$
and let $R_1,..,R_h$ be the level spheres coming from critical
values for $\beta$ between $S$ and $A$. (See Figure \ref{fig:top is
thinnest}).
\begin{figure}
\begin{center} \includegraphics[scale=.3]{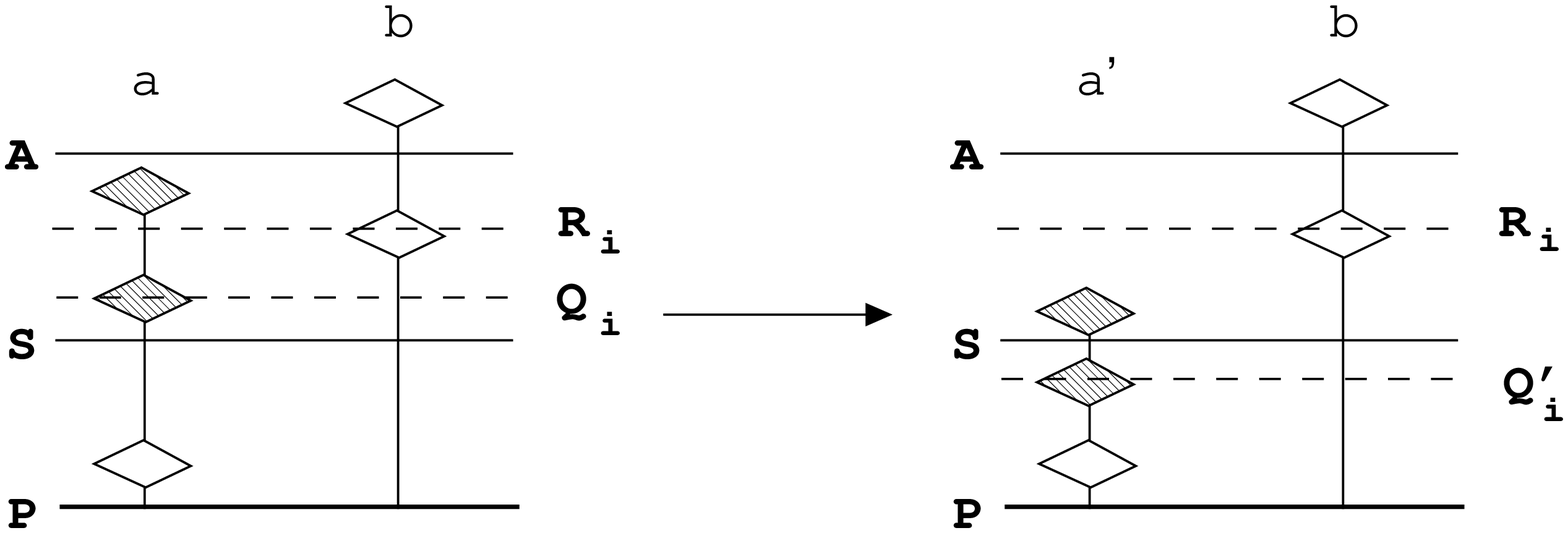} \end{center}
\caption{A rhombus represents a region of $\alpha$ or $\beta$ between two
adjacent thin spheres: a series of minima followed by a series of
maxima.} \label{fig:top is thinnest}
\end{figure}
As
$\beta \cap p^{-1}(s,a)$ was assumed not to be a product, $h \geq 1$.  Push down
the critical points of $\alpha \cap p^{-1}(s, a)$ past $S$. We
will denote by $L'$, $\alpha'$, and $Q_i'$ the images of $L$,
$\alpha$,
and $Q_i$ after the isotopy. We will show that the width of $L'$ is
lower than the width of $L$, leading to a contradiction. The widths
of all level spheres other than $R_i$ and $Q_i$ remain
unchanged by the isotopy. Notice that $|R_i \cap
\alpha'|=0$ as all $R_i$ are above $S$. So we have
$$|R_i \cap L|=|R_i \cap \alpha|+|R_i \cap \beta|>|R_i \cap
\beta|=|R_i \cap \alpha'|+|R_i \cap \beta|=|R_i \cap L'|.$$ The
inequality here is strict as $R_i \cap \alpha \neq \emptyset$. To
show $L'$ is thinner than $L$, we just need to show the
intersection with the $Q_i$ has not increased. For each $Q_i$,
$$|Q_i \cap L|=|Q_i \cap \alpha|+|Q_i \cap \beta|\geq |Q_i \cap
\alpha|+|S\cap \beta|=|Q_i' \cap \alpha'|+|Q_i' \cap \beta|=|Q_i'
\cap L'|.$$ This inequality uses the fact that $S$ provides the
highest counterexample. As $L$ was assumed to be in thin position,
we have a contradiction.
\end{proof}

The following is a direct consequence of Lemma
\ref{lemma:Top alternating sphere is thinnest}.
\begin{corollary}\label{cor:not compressible on one side}
Suppose $P$ is a thin sphere such that for all level spheres $S$
on one side of P  we have that $|P\cap L|\leq |S\cap L|$. Then $P$
is not compressible on that side.
\end{corollary}

Wu's result that the level sphere of minimum width is
incompressible \cite{Wu} is a special case of Corollary
\ref{cor:not compressible on one side}.

By extending the proof of
Lemma \ref{lemma:Top alternating sphere is thinnest} we can
demonstrate that a statement similar to this lemma holds for all
alternating thin level spheres.

\begin{theorem} \label{thm:Alternating spheres are thinner than ones
below} Suppose $C=p^{-1}(c)$ is any alternating thin sphere for
$D$ and $S=p^{-1}(s)$ is any level sphere with $s \in [0,c)$. Then
$|C\cap \alpha|<|S\cap \alpha|$ unless $\alpha \cap\ p^{-1}[s,c]$
is a product, and similarly for $\beta$.
\end{theorem}
\begin{proof}
We know that the theorem holds for the highest alternating sphere
$A$. This follows for $\bbb$ from Lemma \ref{lemma:Top alternating
sphere is thinnest} and for $\alpha$ from the observation that $A\cap
\alpha =\emptyset$. Suppose $C=p^{-1}(c)$ is the highest
alternating sphere for which the statement is not true: i.e. there
exists $S=p^{-1}(s)$ with $s \in [0,c)$ a level sphere such that
either $|C\cap \alpha|\geq |S\cap \alpha|$ and $\alpha \cap\
p^{-1}[s,c]$ is not a product or $|C\cap \beta|\geq |S\cap \beta|$
and $\beta \cap\ p^{-1}[s,c]$ is not a product, say the former.
Take $S$ to be the highest such sphere. As in Lemma \ref{lemma:Top
alternating sphere is thinnest} we may assume that $S$ is thin by
passing to an adjacent thin sphere. By hypothesis, if
$B=p^{-1}(b)$ is the first alternating thin sphere above $C$ then $|B\cap
\beta|<|S\cap \beta|$. $B$ and $C$ are adjacent alternating
spheres, so by Remark \ref{rmk:product region} one of $\alpha \cap
p^{-1}(c,b)$ and $\beta \cap p^{-1}(c,b)$ has to be a product. If
$\alpha \cap p^{-1}(c,b)$ were a product, then $|\alpha \cap
B|=|\alpha \cap C| \geq |\alpha \cap S|$, contradicting the choice of $C$
and $S$. So $\beta \cap p^{-1}(c,b)$ is a
product. Let $Q_1,..,Q_k$ be the level spheres coming from
critical values for $\alpha$ between $S$ and $B$, and let $R_1,..,R_h$
be the level spheres coming from critical values for $\beta$ in the
same region. As
$C$ is alternating, $\beta$ has a critical point directly below
$C$. Thus for all $Q_i$ the region of $\beta$ between $B$ and $Q_i$
is not a product (Figure \ref{fig:Alternating sphere is
thinnest}), and so by the hypothesis $|B\cap \beta|<|Q_i\cap
\beta|$. Push down the critical points of $\beta \cap p^{-1}[s,
b]$ past $S$. Let $\beta'$, $R_i'$
and $L'$ be the images of $\beta$, $R_i$ and $L$ after the
isotopy. The only level spheres whose widths are affected by the
isotopy are the $R_i$ and $Q_i$.

\begin{figure}
\begin{center} \includegraphics[scale=.3]{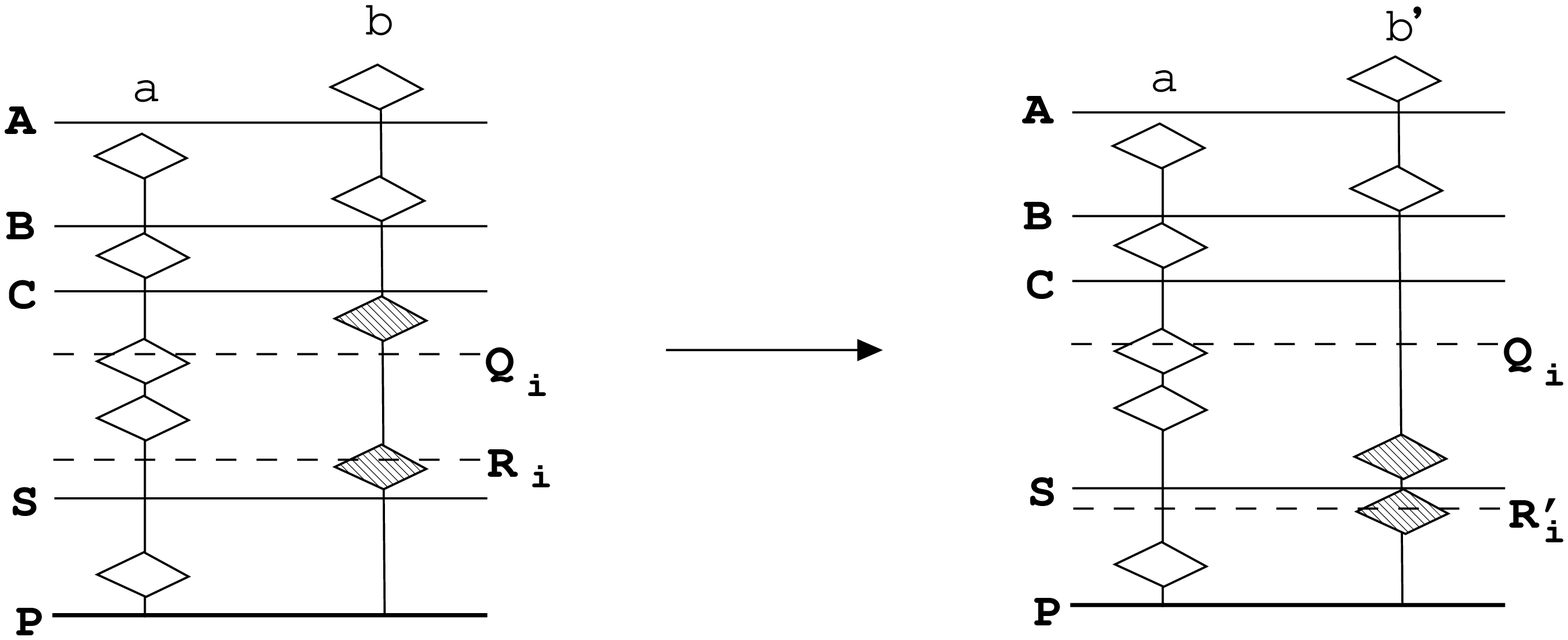} \end{center}
\caption{} \label{fig:Alternating sphere is thinnest}
\end{figure}

Computing the width of $L'$ we obtain, for each $Q_i$ and $R_i$,
$$|Q_i \cap L|=|Q_i \cap \alpha|+|Q_i \cap \beta|> |Q_i \cap \alpha|+|B\cap
\beta|=|Q_i \cap \alpha|+|Q_i\cap \beta'|=|Q_i\cap L'|$$ and
$$|R_i \cap L|=|R_i \cap \alpha|+|R_i \cap \beta|\geq |S \cap
\alpha|+|R_i' \cap \beta'|=|R_i' \cap L'|.$$ To see the second
inequality, note that if for some $R_i$ we have $|R_i \cap
\alpha|<|S \cap \alpha|\leq|C\cap \alpha|$ then a thin
sphere $S'$ adjacent to $R_i$ should have been substituted for $S$ which was assumed
to provide the highest counterexample. Thus the isotopy thins $L$
which was assumed to be in thin position, a contradiction.

\end{proof}

\begin{corollary}\label{cor:alt thin sphere thinner than all below}
If $C=p^{-1}(c)$ is an alternating thin level sphere for $D$, then
its width is lower than the widths of all level spheres lying in $p^{-1}[0,c)$.
\end{corollary}
\begin{corollary}\label{cor:widths monotone decreasing}
The widths of the alternating thin spheres above $P$ are
monotone decreasing.
\end{corollary}

Corollaries \ref{cor:alt thin sphere thinner than all below} and
\ref{cor:widths monotone decreasing} provide restrictions on which
thin spheres can be alternating. As we know that the strands of the
link on one side of a
compressing disk must be a product between adjacent alternating thin
spheres, we can now give a description of how pairs of compressing
disks on opposite sides of $P$ must lie with respect to the link.

\section{Thin level spheres compressible on both sides}

The compressing disk $D$ for $P$ separates $p^{-1}[0, \infty)$
into two balls. Define the \textit{interior ball of D}, denoted
$B^{int(D)}$, to be the closed ball that contains the critical point of $L$
closest to $P$. The \textit{exterior ball of D} will be
denoted $B^{ext(D)}$ and is also closed. Notice that $B^{int(D)}$
and $B^{ext(D)}$ intersect $P$ in disks with boundary $\bdd D$
which we label $D^{int}$ and $D^{ext}$ respectively.

\begin{theorem}\label{thm:Disjoint Interiors}
If $P$ is a thin sphere compressible above and below with
compressing disks $D_u$ and $D_l$ respectively, then $D^{int}_u\cap
D^{int}_l\neq\emptyset$.
\end{theorem}

\begin{proof}
Suppose we have compressing disks  $D_u$ above $P$ and $D_l$ below
$P$ with $D^{int}_u\cap D^{int}_l=\emptyset$. By Remark
\ref{rmk:top is alternating} there is at least one alternating
sphere on each side of $P$. Let $A_u=p^{-1}(a_u)$ be the lowest
alternating thin sphere for $D_u$ and $A_l=p^{-1}(a_l)$ be the
highest alternating thin sphere for $D_l$ (Figure \ref{fig:Disjoint
Interiors}). By Corollary \ref{cor:alt thin sphere thinner than all below},
$$w(P)>w(A_u).$$ By Remark \ref{rmk:product region}, and after
perhaps an $h$-isotopy described in Lemma \ref{lemma:Wu}, $B^{ext(D_u)}\cap p^{-1}[0,a_u]$
is a product, so we can extend $D_l$ up to $A_u$ using this product structure.
Notice that $P$ is now an alternating thin sphere for the extended
$D_l$ and therefore $$w(P)<w(A_u).$$ Thus we have a contradiction.

\begin{figure}
\begin{center} \includegraphics[scale=.3]{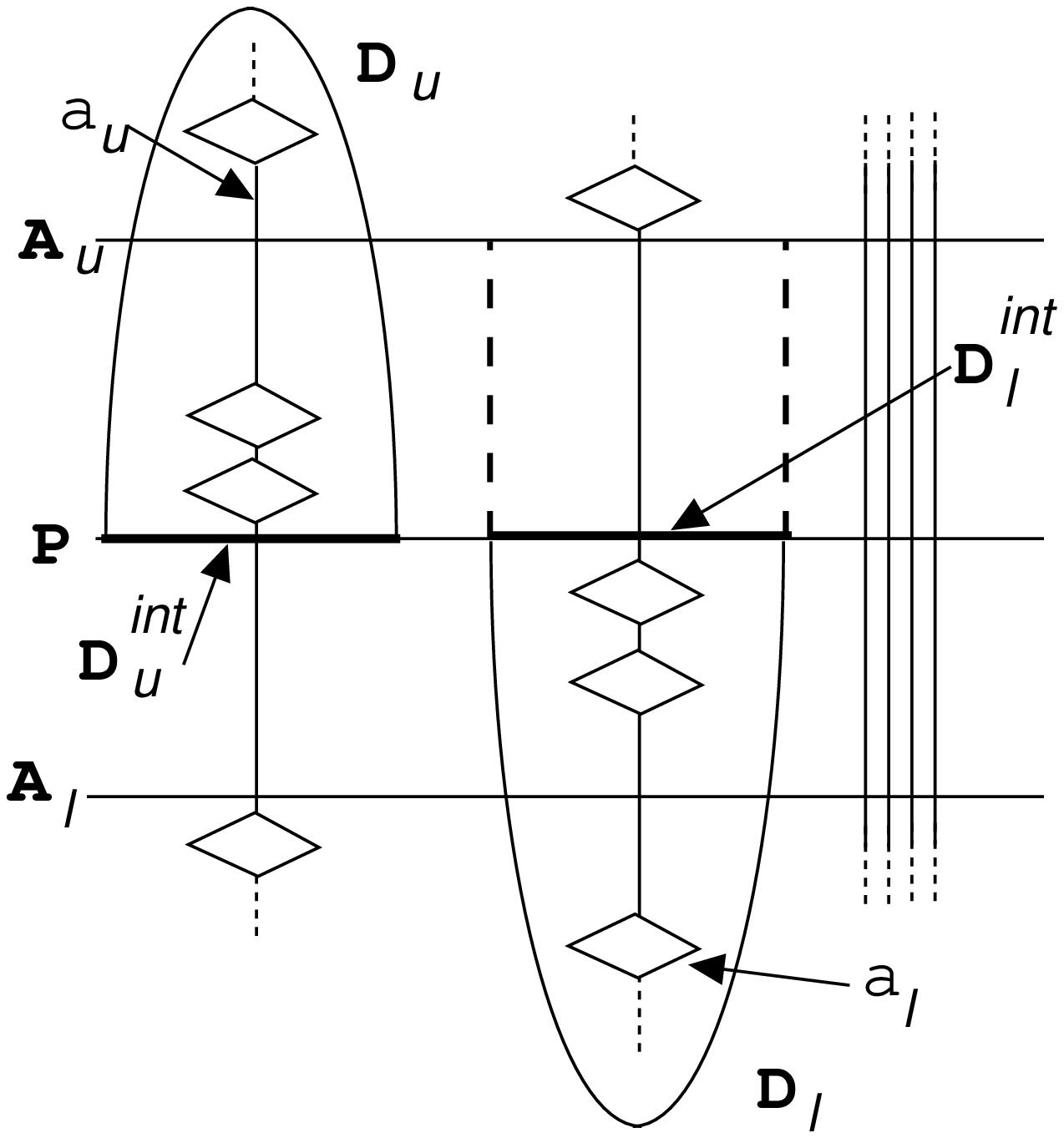} \end{center}
\caption{} \label{fig:Disjoint Interiors}
\end{figure}

\end{proof}

\begin{defin}
A level sphere is strongly compressible in the link complement if
there exist compressing disks $D_u$ and $D_l$ above and below $P$
respectively, such that $\bdd D_u \cap \bdd D_l=\emptyset$. If a
level sphere is not strongly compressible, then it is weakly
incompressible.
\end{defin}

\begin{corollary}\label{cor:nesting}
If a level sphere $P$ is strongly compressible with compressing
disks $D_u$ and $D_l$, then the disks $D^{int}_u$ and $D^{int}_l$
are nested.
\end{corollary}
In fact, using the notation in the proof of Theorem
\ref{thm:Disjoint Interiors}, the relative widths of $A_u$ and
$A_l$ determine the direction of the nesting as the next theorem shows.

\begin{theorem}\label{thm:nesting direction}
$D^{int}_u \subset D^{int}_l$ if and only if $w(A_l)<w(A_u)$.
\end{theorem}
\begin{proof}
If $D^{int}_u \subset D^{int}_l$ then $\bdd D_l$ lies in
$D^{ext}_u$, so just as in the proof of Theorem \ref{thm:Disjoint
Interiors}, we can extend $D_l$ to $A_u$. $A_l$ is still an
alternating thin sphere for the extended disk which is now a
compressing disk for $A_u$. Therefore $w(A_l)<w(A_u)$. The other
direction is obtained by switching the labels $u$ and $l$.
\end{proof}

\begin{corollary}\label{cor:second thinnest}
    If $P$ is a thin sphere with second lowest width among all thin
    spheres, then $P$ is weakly incompressible.
    \end{corollary}

    \begin{proof}
Suppose $P$ is strongly compressible and let $A_u$ and $A_l$ be as in
the theorem.
    From Theorem \ref{thm:nesting direction} we can deduce that
    it is not possible to have $w(A_u)=w(A_l)$. If $P$ is a sphere of second lowest width
    then $A_u$ and $A_l$ will both be of minimal
    width and thus we would have $w(A_u)=w(A_l)$, a contradiction.
\end{proof}

The next theorem provides another sufficient condition for a thin sphere to be weakly
incompressible. First notice the following:

\begin{remark}\label{rmk:compressing disks extend}
\textnormal {Suppose $P$ is a thin sphere compressible above with compressing
disk $D$ and let $Q$ be the first thin sphere above $P$ satisfying
$w(Q)<w(P)$. Then there are no alternating thin spheres for $D$ below
$Q$ and
thus $\alpha$ or $\beta$ must be a product between $P$ and $Q$.
The analogous statement holds for $P$ compressible below.}
\end{remark}

\begin{theorem}\label{thm:weakly incompressible}
Let $P_i$ and $P_k$ be two thin spheres with $P_i$ lying below
$P_k$. Suppose $P_i$ is not compressible above and $P_k$ is not
compressible below. Let $P_j$ be the thin sphere of minimal width
among all spheres strictly between $P_i$ and $P_k$. Then $P_j$ is
weakly incompressible.
\end{theorem}

Notice that in particular the hypothesis holds if
$P_i$ and $P_k$ are both of minimal width.
\begin{proof}
Assume to the contrary that $P_j$ is strongly compressible with
compressing disks $D_u$ and $D_l$ above and below $P_j$
respectively. By Corollary \ref{cor:nesting} the disks $D^{int}_u$
and $D^{int}_l$ are nested. We will assume $D^{int}_u\subset
D^{int}_l$, the other case is similar. By Remark
\ref{rmk:compressing disks extend} we know that $L \cap B^{ext(D_u)}$ must have a product
structure between $P_j$ and $P_k$. As
$\bdd D_l \subset D^{ext}_u$ we can use the product structure to
extend $D_l$ up to $P_k$. But $P_k$ was assumed to be
incompressible from below, so we reach a contradiction.
\end{proof}

\begin{defin}  The region between two level spheres is turbulent if
every strand in the region has at least one critical point.
\end{defin}

\begin{corollary} If the region between two level spheres is turbulent then
a thin sphere of minimal width between them is weakly incompressible.
\end{corollary}

In other words, a thinnest sphere between two other level spheres,
if it is strongly compressible, ``calms" the region between them.

\section{Bound on the number of disjoint irreducible compressing disks}

Once we have placed the link in thin position, the width function
restricted to the thin levels gives a sequence of even integers.
Consider these integers as a set discarding repeating values, arrange
them in (strictly) increasing order and label them $w_0,
w_1,...,w_k$. Thus the lowest width thin sphere has $w_0$ intersection points
with the link,
the  second lowest width thin sphere has $w_1$ intersection points, etc. Of course there may be
several thin spheres all with width $w_i$ and we have no control
in what order the $w_i$ appear when we look at the natural
ordering of the thin spheres given by the height function on
$S^3$.

Using this language Wu's result can be restated as:
\begin{thm}
If $w(P)=w_0$ then $P$ has no compressing disks.
\end{thm}

We can also restate Corollary \ref{cor:second thinnest}:
\begin{cor}
    If $w(P)=w_1$ then $P$ is weakly incompressible.
    \end{cor}

\begin{figure}
 \begin{center} \includegraphics[scale=.3]{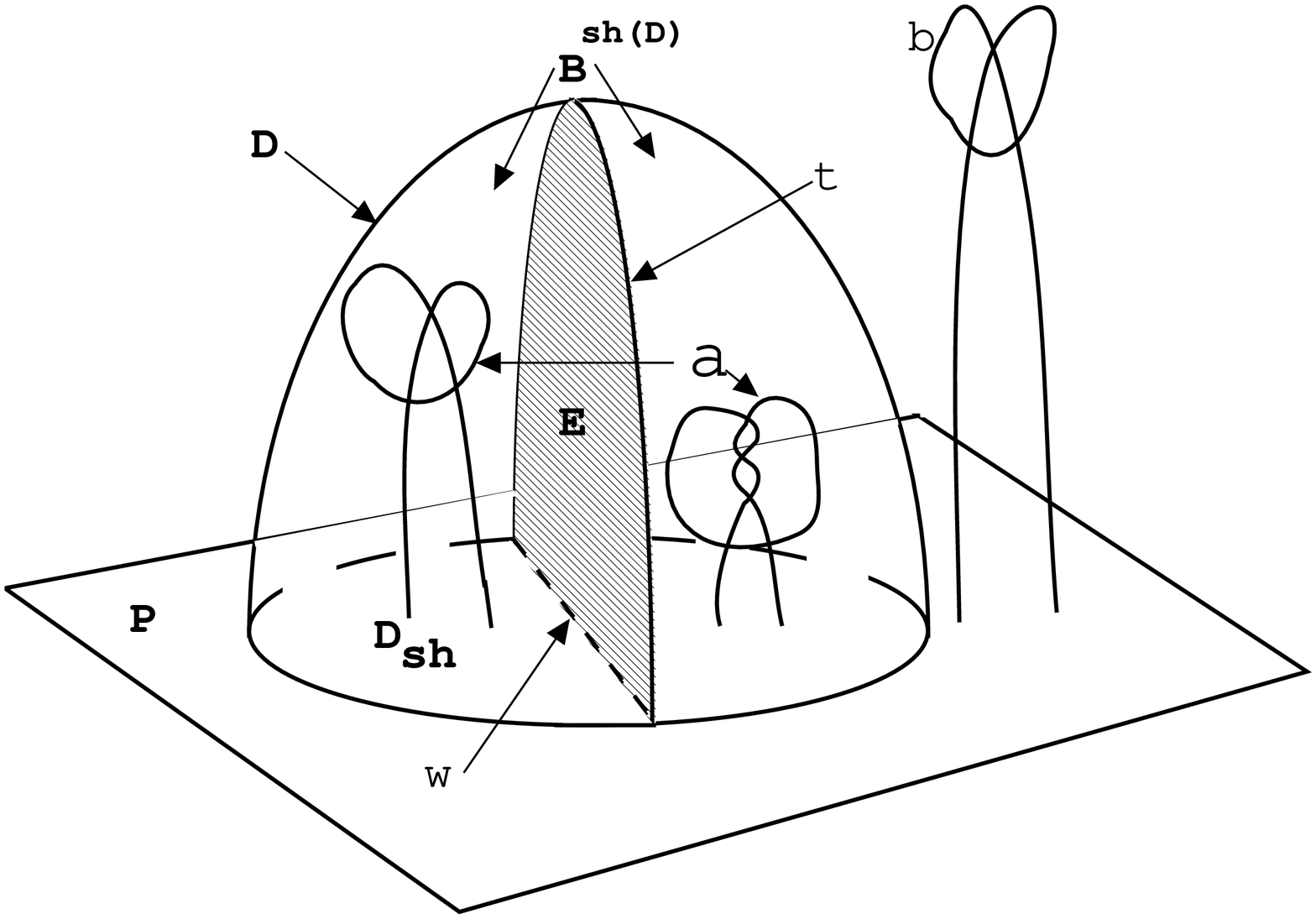} \end{center}
 \caption{} \label{fig:reducingdisk}
 \end{figure}

For the remainder of the paper we will assume that $D$ is a
compressing disk for $P$ lying above, $\alpha$ and $\beta$ are the
strands of a (unsplit) link $L$ contained in each of the two balls cobounded by $D$ and $P$
where $\alpha$ is
shorter than $\beta$. Let the \textit{short ball for $D$},
$B^{sh(D)}$, be the closed ball bounded by $D\cup P$ containing
$\alpha$, and let
$B^{t(D)}$ be the closed ball bounded by $D\cup P$ containing $\beta$.
We will denote the disk $B^{sh(D)}\cap P$ by $D_{sh}$ (See Figure
\ref{fig:reducingdisk}) .

\begin{defin}
A reducing disk for $D$,
$E$, is an embedded disk contained in $B^{sh(D)}$ such that
$\bdd E=\tau \cup \omega$, where $\tau \subset D$ and $\omega
\subset D_{sh}$ is essential in $D_{sh}-L$. $D$ is reducible if such a disk exists and irreducible
otherwise. (Figure \ref{fig:reducingdisk})
\end{defin}

We first prove two straightforward lemmas describing how two compressing disks
can lie relative to each other.

\begin{lemma}\label{lem:noncontainment}
Suppose $D$ and $D'$ are two disjoint, nonparallel,
compressing disks for $P$ such that $D\subset B^{sh( D')}$. Then
$D'$ is reducible.
\end{lemma}
\begin{proof}
Consider an embedded arc $\nu$ spanning the annulus $D'_{sh}-D_{sh}$
and missing the link. A regular neighborhood of $B^{sh(D)}\cup \nu$ 
has boundary a 2-sphere which intersects $B^{sh(D')}$ in a reducing 
disk for $D'$.
\end{proof}

\begin{lemma}\label{lem:no circles}
Suppose $D_1$ and $D_2$ are two compressing disks for $P$. Then there is
an isotopy of $D_2$ in the link complement, fixing $\bdd D_2$,
at the end of which the two disks do not have circles of intersection.
\end{lemma}

\begin{proof}
   Suppose $c \in D_1 \cap D_2$ is a circle
of intersection innermost on $D_1$. The disk on
$D_1$ bounded by $c$, $F_1$, has an interior  disjoint from $D_2$. There
is also a disk $F_2$ in $D_2$ with boundary $c$. Because the link is
not split the embedded sphere $F_1 \cup F_2$
bounds a ball on one side. Therefore $F_2$  can be isotoped
across the ball to remove $c$ from the intersection.

\end{proof}

If $P=p^{-1}(0)$ is a thin sphere we call a thin sphere $A=p^{-1}(a)$
above it
\textit{potentially alternating for $P$} if it satisfies the condition
for the level sphere $C$ in
Corollary \ref{cor:alt thin sphere thinner than all below}, i.e.
its width is strictly lower than the width of any thin sphere in the region
$p^{-1}[0,a)$. Notice that the definition of an alternating sphere
depends on the disk under consideration while the definition of
potentially alternating depends only on the initial sphere $P$. However, as
the term suggests, a sphere can be alternating for some compressing
disk for $P$ only if it is
potentially alternating. We will denote the potentially alternating
thin spheres for $P$ by $A_i$, $i \geq 1$, in order of ascending 
heights (and
descending widths). Of course
if $w(P)=w_k$ then $P$ can have at most $k$ potentially alternating
thin spheres lying above it.

\begin{remark}\label{rmk:product}
\textnormal{If $D$ is a compressing disk for $P$ and $\alpha$ intersects two
  adjacent potentially alternating thin spheres $A_{j-1}$ and $A_j$,
  then one of $\alpha$ or $\beta$ is a product in the region between them.}
    \end{remark}
\begin{defin}
Suppose $D$ is a compressing disk for $P$ so that $\alpha \cap A_{i-1}
\neq \emptyset$ but $\alpha \cap A_i = \emptyset$ then we say that $D$
has height $i$.
\end{defin}
Now we can prove a theorem that has many corollaries including the
main result of this section - a bound on the number of disjoint irreducible compressing disks a thin
level sphere can have.

\begin{theorem}\label{thm:same or disjoint}
  Suppose $D$ and $D'$ are two irreducible compressing disks for
  $P$, and $\alpha$, $\alpha'$ are the strands of $L$ lying in the
  corresponding short balls. Then $height(D)=height(D')$ implies $\alpha=\alpha'$.
  Otherwise $\alpha \cap \alpha'=\emptyset$.
 \end{theorem}

 \begin{proof}
     First we will show that  either $\alpha=\alpha'$ or $\alpha \cap
     \alpha'=\emptyset$. If $D$ and $D'$ are disjoint the result
     follows from Lemma \ref{lem:noncontainment}.

     Let $\Lambda=D\cap D'$
     and assume the disks have been isotoped in the link complement to
     minimize $|\Lambda|$; in particular, by Lemma \ref{lem:no
     circles} $\Lambda$ contains only arcs. By Theorem \ref{lemma:Wu} we can assume
     $D'$ is vertical. Let $\{R_i\}$ be the set of all components of $D-\Lambda$  that lie
     inside $B^{sh(D')}$. Suppose $R_1$ is bounded by curves
     $\lambda_1,\ldots,\lambda_n \in \Lambda$ and
     $\omega_1,\ldots,\omega_n \subset \bdd D\cap D_{sh}'$.
     We will assume $\omega_1$ is outermost of the $\omega_i$ in
     $D_{sh}'$ and thus cobounds a bigon with $\bdd D'$. (See Figure
     \ref{fig:same or disjoint}a.) $R_1$
     separates $B^{sh(D')}$ into closed balls which we will call $F_1$
     and $E_1$ with the labels chosen
      so that the bigon cobounded by $\omega_1$ and $\bdd
     D'$  lies in $F_1$.  This bigon must intersect the link, otherwise we
     could decrease $|\Lambda|$ via an isotopy of $D'$
     pushing $\omega_1$ and any other arcs of $\bdd D \cap D_{sh}'$
     contained in the bigon across it and thus eliminating at least two
      points of $\bdd D \cap \bdd D'$. We can thus conclude that $\alpha'
      \cap F_1 \neq \emptyset$. We will show that in fact
 $\alpha \subset F_1$.

     If we place the point at infinity in
     $D_{sh}' \cap E_1$ we obtain Figure \ref{fig:same or disjoint}b, depicting $\bdd
     B^{sh(D')}-\infty$. The disk $F_1\cap \bdd B^{sh(D')}$ has been shaded.

    \begin{figure}
     \begin{center} \includegraphics[scale=.35]{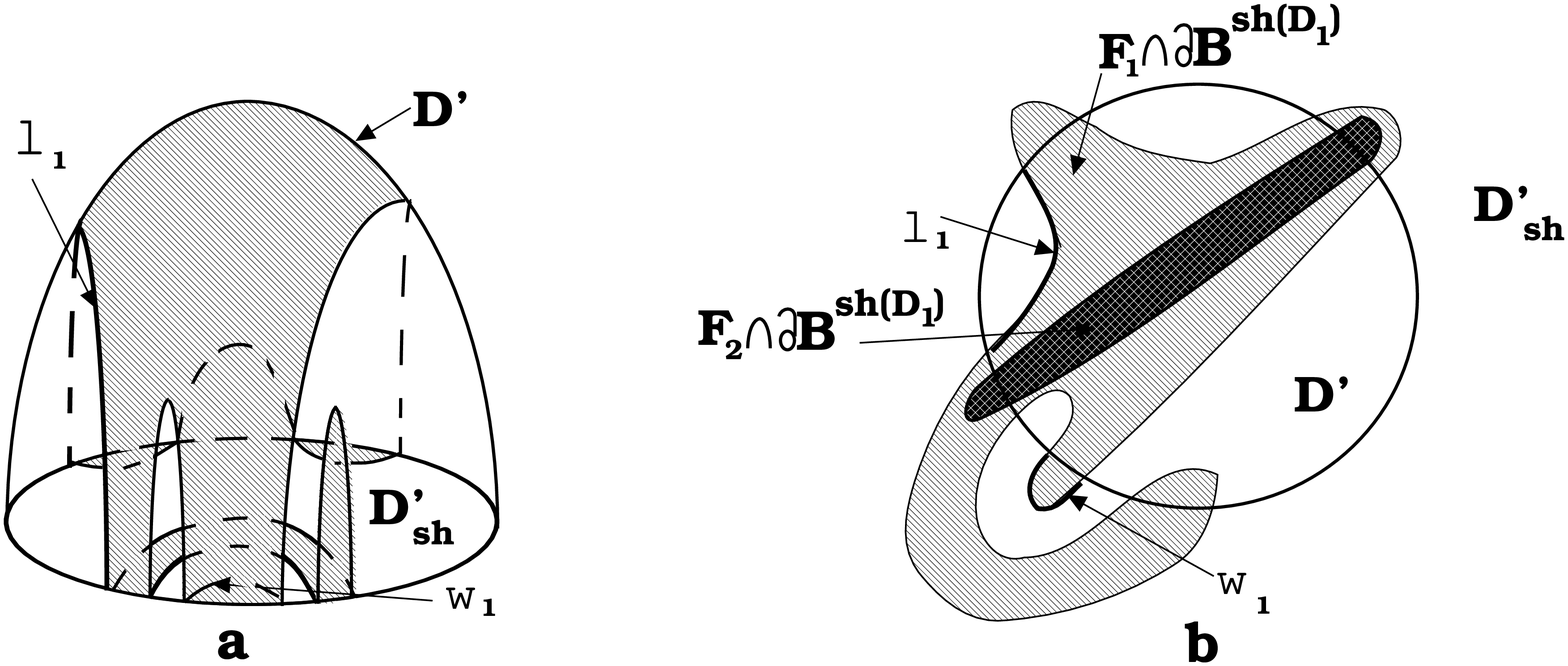} \end{center}
     \caption{} \label{fig:same or disjoint}
     \end{figure}

     Pick $\lambda_1$ to be outermost of the $\lambda_i$ in $D'$ so it cobounds a
     bigon with $\bdd D'$.  Isotope $R_1$ in the link complement
     pushing $\lambda_1$ across the bigon into $D_{sh}'$. No
     restriction is placed on the shading of the bigon. This process decreases
     $|\bdd R_1 \cap \bdd D'|$ by two, so by a series of such moves
     we can isotope $R_1$ so that $\bdd R_1$ and  $\bdd D'$ are disjoint. Notice that $\omega_1$ was never
     isotoped so it must still lie in $D_{sh}'$ and thus all of
     $\bdd R_1$ lies there.  Now $R_1$ is a disk with $\bdd R_1
     \subset D'_{sh}$ and by Lemma
     \ref{lem:noncontainment} it is either not a compressing disk for
     $P$ or it is parallel to $D'$. In other words  $\alpha'$ is
     contained entirely on one side of $R_1$, and as we noticed
     earlier that $\alpha'\cap F_1 \neq \emptyset$, we
     conclude that $\alpha' \subset F_1$.

     Although we do not need this for the rest of the proof, it is
     interesting to observe that as $L \cap D_{sh}' \cap E_1
     =\emptyset$, all outermost arcs
     $\omega_i$ must have bounded shaded bigons of $D_{sh}'$ or we could
     have decreased $|D \cap D'|$ .

     Now consider a second region $R_2$ separating $B^{sh(D')}$ into
     balls $F_2$ and $E_2$ with the labels chosen analogously to the
     labels for $R_1$ and thus $E_2 \cap
     L = \emptyset$. If $R_2$ is contained
     in $F_1$ we will say that $R_2$ is inside $R_1$. Observe that if $R_2$ is inside $R_1$, $F_2$ must be contained in
   $F_1$: an outermost arc of $\bdd R_2 \cap D_{sh}'$ bounds a bigon
   of $D_{sh}'$ contained in $F_2$ (by definition of $F_2$) and also contained in
   $F_1$, which is on the side of $\bdd R_2$ that does not contain $\bdd
   R_1$. (See Figure \ref{fig:same or disjoint}.)
     Thus we can find $R_n$, an innermost region
  with $\alpha'\subset F_n$ and thus
     entirely on one side of $D$. We conclude that  $\alpha' \subset \alpha$ or $\alpha' \subset
     \beta$. Switching the labels of $D$ and $D'$ we also have that  $\alpha\subset
     \alpha'$ or $\alpha \subset
     \beta'$. As $\alpha\cap \beta=\alpha'\cap \beta'=\emptyset$ we must either have
     $\alpha=\alpha'$ or $\alpha \cap \alpha'=\emptyset$ .

     The height of a disk is determined by the strands of $L$
     inside its short ball, therefore if $D$ and
     $D'$ have different heights we must have the case that $\alpha \cap \alpha'=\emptyset$.
     If $height(D)=height(D')=i$, as
     $D'$ was chosen to be vertical, $D' \cap A_{i-1}$ is a single circle
     bounding a subdisk $\tilde{D}$ of $D'$. Let $\tilde{\alpha} = L \cap B^{sh(\tilde{D})}$ and
     $\tilde{\beta} = L \cap
     B^{t(\tilde{D})}$. As $B^{sh(\tilde{D})} \subset B^{sh(D')}$ we must have that
     $\tilde{\alpha} \subset \alpha'$. $\tilde{D}$ is a compressing
     disk for $A_{i-1}$ of height 1 so $\tilde \alpha$ must have
     some critical points. Remark \ref{rmk:product} implies that $\tilde{\beta}$ cannot
  contain any critical points between $A_{i-1}$ and $A_i$ so all strands of
  $\tilde{\beta}$ intersect $A_i$. By the definition of height we must
  have that $\alpha \cap A_i=\emptyset$. We can
  conclude that $\alpha\cap \tilde{\beta} =\emptyset$ so $\tilde{\alpha}
  \subset \alpha$ and therefore $\alpha'\cap\alpha \supset \tilde{\alpha}$.
  Thus we must have the case that $\alpha=\alpha'$.

\end{proof}

\begin{corollary} \label{cor:must intersect}
 Any two distinct irreducible compressing disks for $P$ of the same height must
 intersect.

\end{corollary}

\begin{proof}
 Suppose $D$ and $D'$ are two disks of the same height. By
 Theorem \ref{thm:same or disjoint} we know that
$\alpha=\alpha'$. So if $D \cap D' =\emptyset$,
possibly after switching labels, $D\subset B^{sh(D')}$. Then by
Lemma \ref{lem:noncontainment} we have that $D$ is parallel to $D'$.
\end{proof}

    \begin{corollary}\label{cor:bound}
    Suppose $P$ is a thin sphere and $w(P)=w_n$. Then $P$ has at most
    $n$ disjoint nonparallel irreducible compressing disks.
    \end{corollary}

 \begin{proof}
  Let $\Delta$ be a maximal collection of  disjoint nonparallel irreducible
  compressing disks for $P$.  Corollary \ref{cor:must intersect}
  implies that $\Delta$
  contains at most one disk of each of the $n$ possible heights
  corresponding to the $n$ potentially alternating thin spheres for $P$.
 Therefore $\Delta$ contains at most $n$ disks.

  \end{proof}
 Based on Theorem \ref{thm:same or disjoint}, if there
 is any compressing disk of height $i$ we can, without
 ambiguity,
 introduce the notation $\alpha_i$ and $\beta_i$ for the strands of $L$
 inside and outside the
 short ball of the disk.  By the same
 theorem we also know that $\alpha_i \subset \beta_j$ for every
 $i\neq j$.

 \begin{corollary}
  If there is some $\alpha_k$ containing critical points of $L$ between $A_{j-1}$ and $A_{j}$ with
  $j< k$  then $P$ does not have a compressing disk of height $j$.
\end{corollary}

\begin{proof}

     Suppose a disk of height $j$ exists. Then $\alpha_j \cap A_{j-1} \neq \emptyset$
 but  $\alpha_j \cap A_{j} = \emptyset$, so $\alpha_j$ must have
 critical points between $A_{j-1}$ and $A_j$. By hypothesis  $\alpha_k$ has
 critical points in this region so by Remark \ref{rmk:product} we know
 that $\beta_k$ must be a
 product there. But Theorem \ref{thm:same or disjoint} gives $\alpha_j \subset
 \beta_k$, leading to a contradiction.
\end{proof}

\begin{figure}
    \begin{center} \includegraphics[scale=.35]{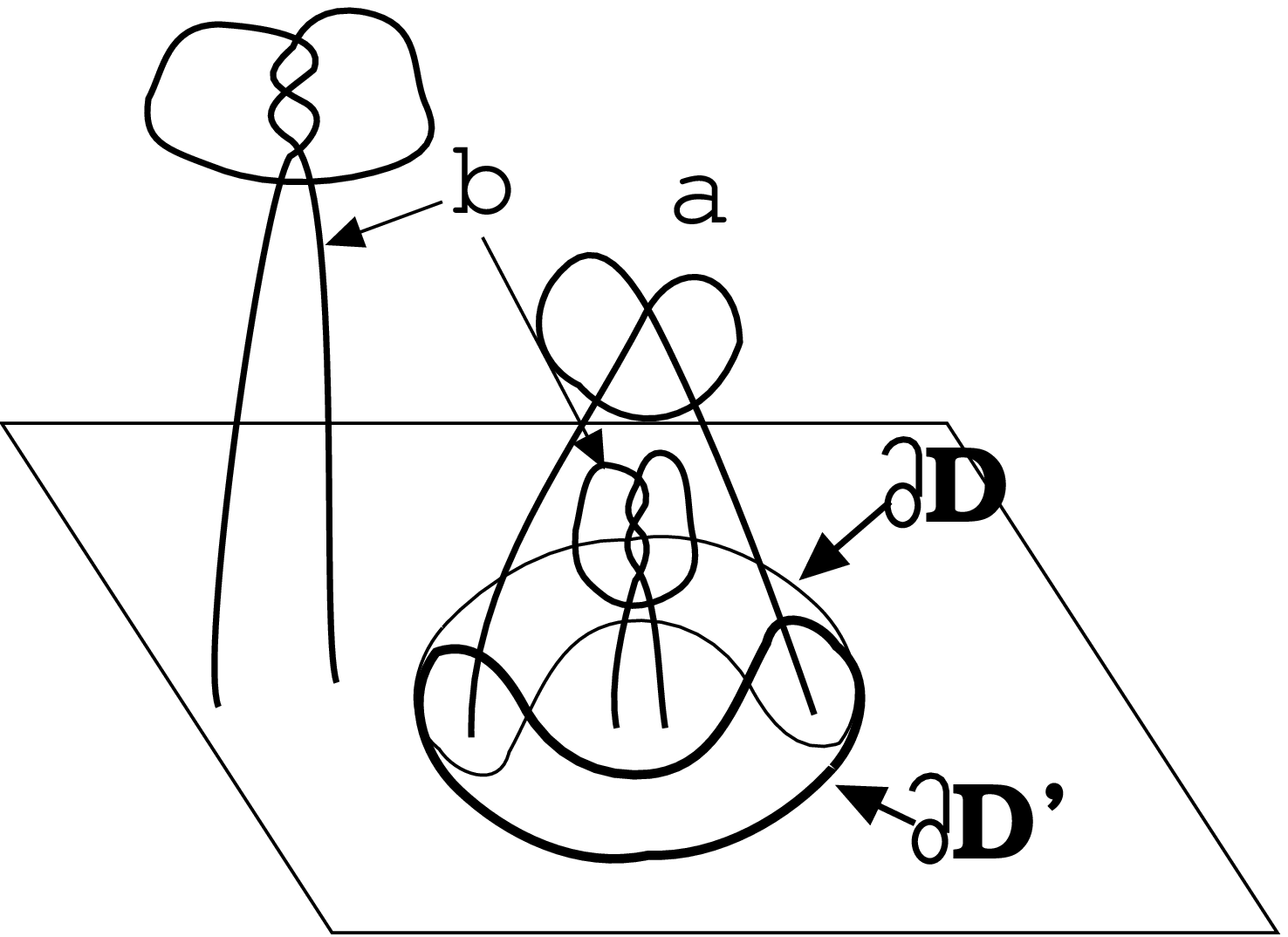} \end{center}
    \caption{} \label{fig:counterexample}
    \end{figure}

Even though the height of the disk completely determines which
strands of $L$ its short ball will contain, it is not true that two
irreducible disks can necessarily be isotoped to be disjoint in the
link complement. Figure \ref{fig:counterexample} shows an example in which
two irreducible disks, $D$ and
$D'$, both of which separate $\alpha$ from $\beta$, cannot be
isotoped to be
disjoint.
However it is possible to pick a family of disjoint irreducible disks
representing all different heights as the next theorem shows.

\begin{theorem} \label{thm:maximal disjoint collection}
 There exists a collection of disjoint irreducible compressing disks
 for $P$ that contains one representative from each possible height.
\end{theorem}

\begin{proof}

Let $\Delta=\cup D_i$ be a maximal in size collection of
    distinct disjoint irreducible compressing disks for $P$.
    By Corollary \ref{cor:must intersect} we know $\Delta$ can contain at most
    one compressing disk of each height.  Suppose there is some
    height $k$ such that $P$ has a compressing disk $D$ of height $k$
    but $\Delta$ does not contain any such disk. Isotope $D$ so that $|D \cap \Delta|$ is minimal. As
    before we will denote by $\alpha$ and $\alpha_i$ the
    intersections of $L$ with the short balls of $D$ and $D_i$
    respectively. As $D \notin
    \Delta$ there exists a curve
 $\lambda \in D\cap \Delta$ outermost in $D$ and thus bounding a subdisk $E$ disjoint from all $D_i$.
 Let $j$ be such that $\lambda
 \in D \cap D_j$ and let $E_j$ be one of the components of $D_j-\lambda$.
 Consider the disk $\tilde D=E \cup E_j$ and let
    $B$ be the ball cobounded by $\tilde D$ and $P$ containing
    $\alpha$.
    Note that $\tilde D \cap \Delta =\emptyset$ so any disk in $\Delta$
    is contained either entirely outside or entirely inside $B$. As
    each $D_i$ is irreducible, $\alpha_i\cap B\neq \emptyset$ implies $\alpha_i
    \subset B$. After renumbering we may assume
    $\alpha_{j_1},\ldots,\alpha_{j_r}$ are contained inside $B$ (this
    collection may be empty). Reduce $\tilde D$
    along a maximal collection of disjoint reducing disks contained in
    $B$ and let $\Delta'$ be the resulting collection of disjoint irreducible disks.
    Note that by Theorem
    \ref{thm:same or disjoint} each $\alpha_{j_i}$, $1 \leq i \leq r$, is
    contained in the short ball of a unique disk in $\Delta'$. Let
    $D'_{j_i}$ be the disk in $\Delta'$ with short ball containing
    $\alpha_{j_i}$. Because $B$ was chosen so that $\alpha \subset B$
    there must be an irreducible disk $D'_{j_0}$ inside $B$
  with a short ball containing $\alpha$. Then the set
    $\{D'_{j_0},\ldots,D'_{j_r}\} \cup (\Delta-\{D_1,\ldots,D_r\})$ contains
   $|\Delta|+1$ disjoint irreducible compressing disks, thus contradicting
   the maximality of $\Delta$.

 \end{proof}

\section{Thin spheres with unique compressing disks}
 In the previous section we restricted our attention to irreducible
 disks. However when discussing disks of height 1 we can drop this
 restriction.
 \begin{lemma}\label{lem:height one always irreducible}
     Any compressing disk of height 1 is irreducible.
  \end{lemma}

  \begin{proof}
      Suppose now $P$ has a compressing disk $D$ of height 1 that is reducible. If
      $\cup_{i=1}^n  E_i$ is a maximal collection of disjoint reducing disks for
      $D$, then the boundary of a regular neighborhood of $D \cup
      (\cup_{i=1}^n E_i)$ contains a collection of $n+1$
      compressing disks of height 1 with disjoint short balls. This contradicts
Theorem \ref{thm:same or disjoint}.
  \end{proof}

\begin{theorem}\label{thm:unique height one disk}
If $P$ has only one potentially alternating thin sphere on some side,
then $P$ has at most one compressing disk on that side.
\end{theorem}

\begin{proof}

If there are counterexamples to this theorem, choose $D_1$ and $D_2$ to
be the pair which, among all counterexamples, can be isotoped to have the fewest 
intersection arcs.
By the hypothesis any compressing disk for $P$ lying above $P$
must necessarily have height 1
so $height(D_1)=height(D_2)=1$.
By Lemma \ref{lem:height
one always irreducible}, $D_1$ and $D_2$ must be irreducible and by
Corollary \ref{cor:must intersect}, they intersect.
Isotope $D_1$ and $D_2$ so as to minimize their intersection, in particular all circles of
intersection are removed by Lemma \ref{lem:no circles}. Let $\lambda$ be an arc of intersection
outermost on $D_1$. Then $\lambda$ bounds a disk $E_1$ in $D_1$ whose interior is disjoint from
$D_2$. Let $E_2$ be one of the components of $D_2-\lambda$.
By minimality of $|D_1 \cap D_2|$
the disk $E=E_1 \cup E_2$ must be a compressing disk for $P$,
necessarily of height 1.  $E$ is disjoint from $D_2$ and has at
least one fewer arc of intersection with $D_1$ than $D_2$ did.
Therefore neither of the pairs $E, D_1$ nor $E, D_2$ are
counterexamples to the theorem at hand, and so $E$ is parallel to both
$D_1$ and $D_2$. Thus $D_1$ and $D_2$ are
parallel to each other and were also not a counterexample to the
theorem.
\end{proof}

\begin {corollary}
If $w(P)=w_1$ then $P$ has at most one compressing disk on each side.
\end{corollary}

\begin{proof}
As $w(P)=w_1$, $P$ can have at most one potentially alternating thin
sphere on each side. The result then follows from
Theorem \ref{thm:unique height one disk}.

\end{proof}
 It is straightforward to observe the following sufficient condition
 for $P$ to not have a compressing disk of height 1.
 \begin{observation}
     If there is a strand of $L$ that has a critical point
     between $P$ and $A_1$ and intersects $A_1$, then $P$ does not have a
     compressing disk of height 1 above.
 \end{observation}

\section{Application to Knots}
\begin{theorem}
Suppose K is a knot or a prime link in thin position and suppose $w_1=w_0+2$. If
$w(P)=w_1$, then $P$ is incompressible.
\end{theorem}
\begin{proof}
Suppose there is a compressing disk $D$ for $P$ lying above it. $P$ has at most one
potentially alternating thin sphere, $A=A_1$, with $w(A)=w_0$ thus
$D$ must have height 1. For a disk
of height 1 we always have that $B^{sh(D)}=B^{int(D)}$ because $L$ can only have critical
points between $P$ and $A_1$ on one side of the disk and that is necessarily the
short ball. As before let
$\alpha$ denote the strands of $K$ contained in $B^{sh(D)}=B^{int(D)}$ and
let $\beta$ denote those contained in $B^{ext(D)}$. Also let $D^{ext}$ and
$D^{int}$ be the two disks that $\bdd D$ bounds on $P$. (See Figure
\ref{fig:Knotsum}.) By the product structure of $\beta$ between
$P$ and $A$ we have that $|\beta \cap A|=|\beta \cap
P|$. As $|\alpha \cap A|=0$ we must have that $|\alpha \cap P|=2$,
and thus the sphere $D \cup D^{int}$ gives a decomposition of $K$ as a
direct sum $K_1\#K_2$ so $K$ cannot be a prime link. Suppose then
that $K$ is a knot and let
$K_1$ be the summand contained in $B^{int(D)}$,
\begin{figure}
\begin{center} \includegraphics[scale=0.35]{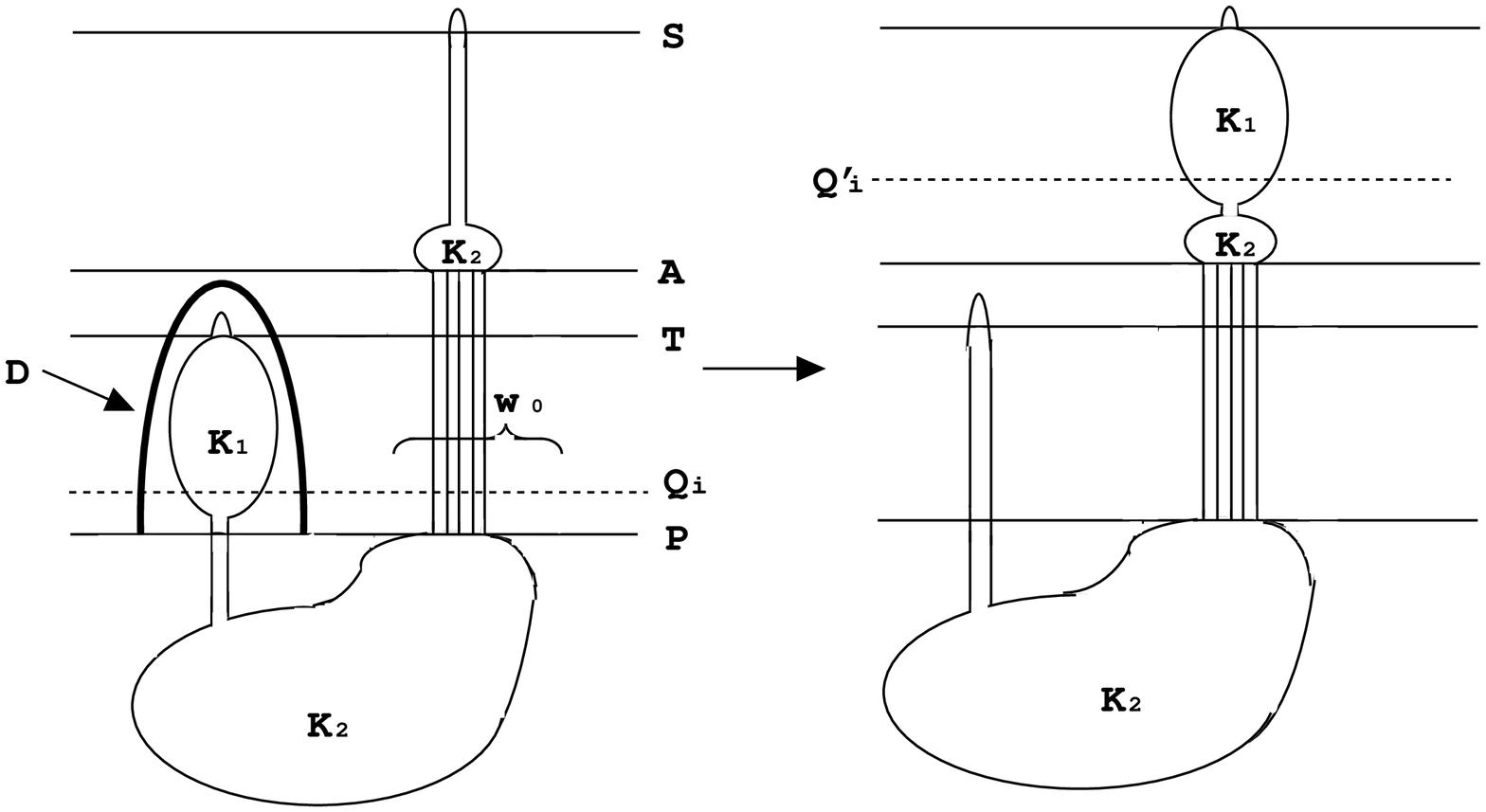} \end{center}
\caption{} \label{fig:Knotsum}
\end{figure}
let $T$ be the
level sphere directly below the highest maximum for $\alpha$ and let
$S$ be the level sphere directly below the highest maximum for
$\beta$. Imagine cutting the region of $\alpha$ between $P$ and
$T$ and inserting it at $S$. As $K$ is a direct sum the same
effect can be achieved via an isotopy making $K_1$ small and
sliding it along $K_2$ to the desired position. No new critical
points have been introduced and the only level spheres affected
are those between $P$ and $T$ which we denote by $Q_i$ before the
isotopy and $Q_i'$ after the isotopy. The knot after the
isotopy will be denoted by $K'$. For all $i$ we have that
$$|Q_i\cap K|=|Q_i\cap \alpha|+|Q_i\cap \beta|>|Q_i\cap \alpha|=|Q_i'\cap
K'|$$ so the width of $K$ has been decreased thus contradicting the
assumption.
\end{proof}

\section*{Acknowledgement}
I would like to thank Marty Scharlemann for many helpful
conversations.
 \nocite{*}
\bibliography{mybiblio}
\bibliographystyle{plain}

\end{document}